\documentclass[12pt,reqno]{amsart}%{article}
\usepackage{amsfonts,amsmath,amssymb}
\usepackage{underscore}

\allowdisplaybreaks
\advance\textheight by \topskip

\usepackage{mathrsfs}

\usepackage[latin1]{inputenc}
\usepackage{graphicx}

\usepackage{tikz}

\usepackage{euscript}

\def\[{[\! [}
\def\]{]\! ]}

\begin{document}

\title{Kn\"odel walks in a B\"ohm-Hornik environment}

\author[H.~Prodinger]{Helmut Prodinger}

\address{Helmut Prodinger,
Mathematics Department, Stellenbosch University,
7602 Stellenbosch, South Africa, and NITheCS (National Institute for
Theoretical and Computational Sciences), South Africa.}
\email{hproding@sun.ac.za}

\date{\today}

\begin{abstract}
Ideas of Kn\"odel and B\"ohm-Hornik about walks in certain graphs, resembling the classical symmetric random walk on the integers,
are combined. All the relevant generating functions (although occasionally quite involved) are made fully explicit.
\end{abstract}


\maketitle

\section{Introduction}

The standard random walk on the non-negative integers may be visualized by the following graph (only the first 8 states are shown):
\begin{figure}[h]

	\begin{center}
		\begin{tikzpicture}[scale=1.5]

			\foreach \x in {0,1,2,3,4,5,6,7,8}
			{
				\draw (\x,0) circle (0.05cm);
				\fill (\x,0) circle (0.05cm);
				%	\node[circle,draw ] (c) at (\x,0){$\x$};
			}

			\fill (0,0) circle (0.08cm);

			\foreach \x in {0,2,4,6}
			{
				\draw[thick,  -latex] (\x,0) to[out=20,in=160]  (\x+1,0);	
				\draw[thick,  -latex] (\x+1,0) to[out=200,in=-20]  (\x,0);	
			}
			\foreach \x in {1,3,5,7}
			{
				\draw[thick,  -latex] (\x,0) to[out=20,in=160]  (\x+1,0);	
				\draw[thick,  -latex] (\x+1,0) to[out=200,in=-20]  (\x,0);	
			}

			\foreach \x in {0,1,2,3,4,5,6,7}
			{
				%	\draw[thick, -latex] (\x+1,0) to  (\x,0);	
				\node at  (\x+0.1,0.15){\tiny$\x$};
			}

			\node at  (8+0.1,0.15){\tiny$8$};

		\end{tikzpicture}
	\end{center}
\caption{Standard symmetric random walk on the non-negative integers}
\end{figure}
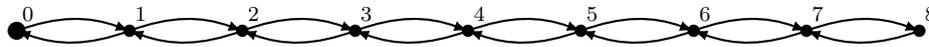

One starts in state 0 and can go up/down one step, each with the same probability.

B\"ohm and Hornik \cite{BH} introduces a related model: up-steps occur with probability $\alpha$ and down-steps occur with probability $\beta=1-\alpha$,
but after each step $\alpha$ and $\beta$ change their roles. The follow graph is useful to grasp the idea. 
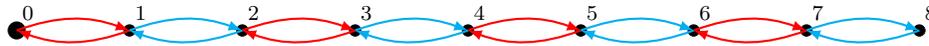
\begin{figure}[h]

	\begin{center}
		\begin{tikzpicture}[scale=1.5]

			\foreach \x in {0,1,2,3,4,5,6,7,8}
			{
				\draw (\x,0) circle (0.05cm);
				\fill (\x,0) circle (0.05cm);
				%	\node[circle,draw ] (c) at (\x,0){$\x$};
			}

			%\draw[thick, red, -latex]    (0,-1) to (-1,0.5);
			
			\fill (0,0) circle (0.08cm);

			\foreach \x in {0,2,4,6}
			{
				\draw[thick, red, -latex] (\x,0) to[out=20,in=160]  (\x+1,0);	
				\draw[thick, red, -latex] (\x+1,0) to[out=200,in=-20]  (\x,0);	
			}
			\foreach \x in {1,3,5,7}
			{
				\draw[thick, cyan, -latex] (\x,0) to[out=20,in=160]  (\x+1,0);	
				\draw[thick, cyan, -latex] (\x+1,0) to[out=200,in=-20]  (\x,0);	
			}

			\foreach \x in {0,1,2,3,4,5,6,7}
			{
				%	\draw[thick, -latex] (\x+1,0) to  (\x,0);	
				\node at  (\x+0.1,0.15){\tiny$\x$};
			}

			\node at  (8+0.1,0.15){\tiny$8$};

		\end{tikzpicture}
	\end{center}
\caption{Red edges are labelled with the weight $\alpha$, blue edges with $\beta$}
\end{figure}

B\"ohm and Hornik \cite{BH} consider  random walks to the non-negative integers and  on
the full set of integers as well.
Alternative/additional analysis can be found in \cite{PP}.

Another twist of a random walk occurs in a model introduced by Kn\"odel \cite{Kn}:
There are bins of size 1 and small items (size $\frac13$) and large items (size $\frac23$) 
arrive with the same probability. States correspond to boxes filled with just one large item each.
There is one exception, when a small item arrives at the origin. In this case, it cannot be used
to complete a partially filled bin, and an extra state is introduced. See \cite{kernel} and some
referenced papers for analysis.

It is the purpose of this paper to combine the ideas of Kn\"odel and B\"ohm-Hornik:
Large items arrive with probability $\alpha$ and small items with probability $\beta$, but after
each step the roles of $\alpha$ and $\beta$ are changed. The graph  with two layers of states
will explain the scenario readily.
\begin{figure}[h]

	\begin{center}
		\begin{tikzpicture}[scale=1.5]
			%\draw (0,0) circle (0.1cm);

			\foreach \x in {0,1,2,3,4,5,6,7,8}
			{
				\draw (\x,0) circle (0.05cm);
				\fill (\x,0) circle (0.05cm);
				%	\node[circle,draw ] (c) at (\x,0){$\x$};
			}
			
			\foreach \x in {0,1,2,3,4,5,6,7,8}
			{
				\draw (\x,-1) circle (0.05cm);
				\fill (\x,-1) circle (0.05cm);
				%	\node[circle,draw ] (c) at (\x,0){$\x$};
			}

			\draw (-1,0.5) circle (0.05cm);	\fill (-1,0.5) circle (0.05cm);
			\draw (-1,-1.5) circle (0.05cm);	\fill (-1,-1.5) circle (0.05cm);

			\draw[thick, cyan, -latex] (0,0) to [out=150,in=-30] (-1,0.5);	
			\draw[thick, cyan, -latex]  (-1,0.5)to [out=00,in=120](0,0);

			\draw[thick, red, -latex] (0,-1) to  (-1,-1.5);	
			\draw[thick, red, -latex]  (-1,-1.5)to [out=0,in=-120](0,-1);

			%\draw[thick, red, -latex]    (0,-1) to (-1,0.5);
			\draw[thick, cyan, -latex]     (-1,-1.5)to (1,0) ;
			\draw[thick, red, -latex]     (-1,0.5)to (1,-1) ;
			
			\fill (0,0) circle (0.08cm);

			\foreach \x in {0,2,4,6}
			{
				\draw[thick, red, -latex] (\x,0) to[out=20,in=160]  (\x+1,0);	
				\draw[thick, red, -latex] (\x+1,0) to[out=200,in=-20]  (\x,0);	
			}
			\foreach \x in {1,3,5,7}
			{
				\draw[thick, cyan, -latex] (\x,0) to[out=20,in=160]  (\x+1,0);	
				\draw[thick, cyan, -latex] (\x+1,0) to[out=200,in=-20]  (\x,0);	
			}

			\foreach \x in {0,2,4,6}
			{
				\draw[thick, cyan, -latex] (\x,-1) to[out=20,in=160]  (\x+1,-1);	
				\draw[thick, cyan, -latex] (\x+1,-1) to[out=200,in=-20]  (\x,-1);	
			}
			
			\foreach \x in {1,3,5,7}
			{
				\draw[thick, red, -latex] (\x,-1) to[out=20,in=160]  (\x+1,-1);	
				\draw[thick, red, -latex] (\x+1,-1) to[out=200,in=-20]  (\x,-1);	
			}
			
\foreach \x in {0,1,2,3,4,5,6,7}
{
	%	\draw[thick, -latex] (\x+1,0) to  (\x,0);	
	\node at  (\x+0.1,0.15){\tiny$\x$};
}

			\node at  (8+0.1,0.15){\tiny$8$};

		\end{tikzpicture}
	\end{center}
\caption{The Kn\"odel-B\"ohm-Hornik graph}
\label{KBH}
\end{figure}
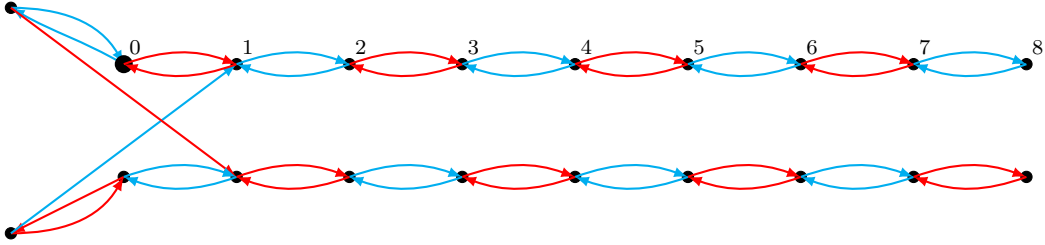
The rest of the paper is devoted to derive generating functions for walks starting at the origin and ending in a prescribed state.
The kernel method \cite{kernel} and the heavy use of computer algebra (Maple) will be essential.

First, we start with a direct approach, which is a brute-force procedure. It leads to four equations, and eventually to biquadratic equations.
Computers are capable of handling this, but the next section mostly serves as an invitation to a more sophisticated approach,
using only two functions (not four). And, lo and behold, after a certain substitution, the ugly beast turns into a beautiful swan.

\section{Brute-force Analysis}

We introduce the following generating functions: $f_i=f_i(z)$ has as coefficient of $z^n$ the probability to reach state $i$
from the upper layer in $n$ steps, starting from the origin (state 0). The function  $g_i$ is similar, but refers to the lower layer of
states. Finally, the extra states and their generating functions are called $P$ resp.\ $Q$.

From the diagram, considering the last step made, one can see the recursions\footnote{One referee suggests to use the symbolic method, as described in Analytic Combinatorics, by Flajolet and Sedgewick. I have known Philippe Flajolet for more than 30 years, being a close co-author as well, and I am quite confident that he would have chosen the kernel method as well, as he often did. It would be a very artificial enterprise to write out symbolic equations. If such an approach would be successful, there would be a symbolic expression for each ``state'', and they would depend on each other recursively.}
\begin{align*}
f_i&=\beta zf_{i-1}+\alpha z f_{i+1}, \ i=2,4,6,\dots,\\
f_i&=\alpha zf_{i-1}+\beta z f_{i+1}, \ i=3,5,7,\dots,\\
f_1&=\alpha zf_{0}+\beta z f_{2}+\beta z Q=\alpha zf_{0}+\beta z f_{2}+\beta\alpha z^2 g_0,\\
f_0&=1+\beta z P+\alpha z f_1=1+\beta^2 z^2 f_0+\alpha z f_1,\\
P&=\beta z f_0,\\
g_i&=\alpha zg_{i-1}+\beta z g_{i+1}, \ i=2,4,6,\dots,\\
g_i&=\beta zg_{i-1}+\alpha z g_{i+1}, \ i=3,5,7,\dots,\\
g_1&=\beta zg_{0}+\alpha z g_{2}+\alpha z P=\beta zg_{0}+\alpha z g_{2}+\alpha\beta z^2 f_0,\\
g_0&=\alpha z Q+\beta z g_1=\alpha^2 z^2 g_0+\beta z g_1,\\
Q&=\alpha z g_0.
\end{align*}
In order to attack this system, we introduce a second variable $u$ and consider the following four bivariate generating
functions:
\begin{align*}
F_e(u)&=\sum_{i\ge0}u^{2i}f_{2i},\quad F_o(u)=\sum_{i\ge0}u^{2i+1}f_{2i+1},\\
G_e(u)&=\sum_{i\ge0}u^{2i}g_{2i},\quad G_o(u)=\sum_{i\ge0}u^{2i+1}g_{2i+1};
\end{align*}
`e' stands for `even', `o' stands for odd. Summing the first recursion, we find (omitting the variable $u$ for the moment)
\begin{equation*}
F_e-f_0=\beta zuF_o+\frac{\alpha z}{u}(F_o-uf_1);
\end{equation*}
adding the recursion for $f_0$ leads to
\begin{equation*}
	F_e=\beta zuF_o+\frac{\alpha z}{u}F_o+1+\beta^2 z^2 f_0.
\end{equation*}
Similarly, for the odd indices
\begin{equation*}
F_o-uf_1=\alpha zu(F_e-f_0)+\frac{\beta z}{u}(F_e-f_0-u^2f_2)
\end{equation*}
and further
\begin{equation*}
	F_o=\alpha zuF_e+\frac{\beta z}{u}(F_e-f_0)+u\beta\alpha z^2 g_0.
\end{equation*}
The same procedure is done for the even indices and the $g_i$'s:
\begin{equation*}
G_e-g_0=\alpha zu G_o+\frac{\beta z}{u}(G_o-ug_1)
\end{equation*}
and
\begin{equation*}
	G_e=\alpha zu G_o+\frac{\beta z}{u}G_o+\alpha^2 z^2 g_0.
\end{equation*}
Finally, for the odd indices
\begin{equation*}
G_o-ug_1=\beta zu(G_e-g_0)+\frac{\alpha z}{u}(G_e-g_0-u^2g_2)
\end{equation*}
and
\begin{equation*}
	G_o=\beta zuG_e+\frac{\alpha z}{u}(G_e-g_0)+u\alpha\beta z^2 f_0.
\end{equation*}
For the reader's convenience we collected the four equations that we (and Maple) have to deal with:
\begin{align*}
	F_e&=\beta zuF_o+\frac{\alpha z}{u}F_o+1+\beta^2 z^2 f_0,\\
	F_o&=\alpha zuF_e+\frac{\beta z}{u}(F_e-f_0)+u\beta\alpha z^2 g_0,\\
	G_e&=\alpha zu G_o+\frac{\beta z}{u}G_o+\alpha^2 z^2 g_0,\\
	G_o&=\beta zuG_e+\frac{\alpha z}{u}(G_e-g_0)+u\alpha\beta z^2 f_0;
\end{align*}
we note again that $f_0=F_e(0)$ and $g_0=G_e(0)$.

Maple can solve this, but the solution is implicit since it still depends on $f_0$ and $g_0$. The expressions are quite long, and
they all share the same denominator $D$:
\begin{equation*}
D={u}^{2}-{z}^{2}{u}^{4}\alpha-{z}^{2}{u}^{2}+2 {u}^{2}\alpha{z}^{2}+{\alpha}^{2}{z}^{
	2}{u}^{4}-2 {u}^{2}{\alpha}^{2}{z}^{2}-{z}^{2}\alpha+{z}^{2}{\alpha}^{2}.
\end{equation*}

Then
$DF_e={u}^{2}- \alpha f_0 {z}^{2}+{\alpha}^{2}f_0 {z}^{2}-2 {z}^{3}{u}^{4}{
	\alpha}^{2}g_0+{z}^{3}{u}^{4}\alpha g_0+{z}^{3}{u}^{2}{\alpha}^{2}g_0+{
	z}^{3}{\alpha}^{3}{u}^{4}g_0-{z}^{3}{u}^{2}{\alpha}^{3}g_0
$
and
$
DF_o=uz ( -{u}^{2}{\alpha}^{2}zg_0-f_0+\alpha f_0+{u}^{2}\alpha zg_0-{z}^{2}{\alpha}^{3}f_0+\alpha{u}^{2}-\alpha+{\alpha}^{3}{z}^{2}{u}^{2}f_0+1+{u}^{2}\alpha{z}^{2}f_0-2 {u}^{2}{\alpha}^{2}{z}^{2}f_0+f_0 {z}^{	2}-3 \alpha f_0 {z}^{2}+3 {\alpha}^{2}f_0 {z}^{2} )
$
and
$DG_e=-\alpha{z}^{2} ( -\alpha{u}^{4}zf_0+{\alpha}^{2}{u}^{4}zf_0+g_0-{
	u}^{2}zf_0+2 {u}^{2}\alpha zf_0-\alpha g_0-{u}^{2}{\alpha}^{2}zf_0
)
$
and finally
$DG_o=
-uz\alpha ( -{u}^{2}\alpha{z}^{2}g_0+g_0-{u}^{2}zf_0+{u}^{2}
\alpha zf_0+{u}^{2}{\alpha}^{2}{z}^{2}g_0-{\alpha}^{2}{z}^{2}g_0
)
$.

The denominator $D$ has 4 roots, considering $u$ as the variable:
\begin{equation*}
s_1={\frac {\sqrt {\alpha ( 1-\alpha )  ( 1-2 {z}^{2}{				\alpha}^{2}+2 {z}^{2}\alpha-{z}^{2}-\sqrt { ( 1-z )  ( 1+z)  ( 1-z+2 z\alpha )  (1+z-2 z\alpha ) }
			) }}{\sqrt {2}z\alpha ( 1-\alpha ) }},
\end{equation*}
\begin{equation*}
s_2=-s_1,\ s_3=\frac1{s_1},\ s_4=\frac1{s_2}.
\end{equation*}

The factors $u-s_1$ and $u-s_2$ are `bad' in the sense of the kernel method \cite{kernel}, i.~e., they don't lead to a power series expansion around the origin.
Consequently, the numerators of the four functions must be divisible by both factors. Applying this principle to $F_e$ and $G_e$ leads to two equations, from
which $f_0$ and $g_0$ can be computed. Again, the expressions are long, and an auxiliary quantity $W$ is used:
\begin{equation*}
W=\sqrt {  (1-z  )   ( z+1  )   (1 -z+2 z\alpha	 )   (1 +z-2 z\alpha ) }.
\end{equation*}
Here are the results:
\begin{equation*}
f_0={\frac {  \Xi_1}{4{
			\alpha}^{2}{z}^{4}  ( -1+z  )   ( z+1  )   ( -1+\alpha
		 ) ^{2}  ( -1+{z}^{2}-3 {z}^{2}\alpha+3 {z}^{2}{\alpha}^{2}  ) 
}}
\end{equation*}
with $\Xi_1=( -3 {z}^{4}{\alpha}^{2}+{\alpha}^{2}W{z}^{2}+3 {z}^{2}{\alpha}^
{2}+3 {z}^{4}\alpha-aW{z}^{2}-3 {z}^{2}\alpha-{z}^{4}+W{z}^{2}+2 {z}^{2}-1-W
)   ( 2 {z}^{2}{\alpha}^{2}-2 {z}^{2}\alpha+{z}^{2}-1+W  ) $\\ and
\begin{equation*}
g_0={\frac {\Xi_2}{8{z}^{7}  ( -1+z  ) 
		 ( z+1  )   ( -1+\alpha  ) ^{4}  ( -1+{z}^{2}-3 {z}
		^{2}\alpha+3 {z}^{2}{\alpha}^{2}  ) {\alpha}^{3}}}
\end{equation*}
with $\Xi_2=  ( -3 {z}^{4}{\alpha}^{2}+{\alpha}^{2}W{z}^{2}+3 {z}^{2}{\alpha}^
{2}+3 {z}^{4}\alpha-aW{z}^{2}-3 {z}^{2}\alpha-{z}^{4}+W{z}^{2}+2 {z}^{2}-1-W
)   ( -{z}^{4}{\alpha}^{2}+2 {z}^{4}{\alpha}^{3}+1-3 {z}^{2}{\alpha}^{2}
-W+{\alpha}^{2}W{z}^{2}-{z}^{2}+2 {z}^{2}\alpha  )   ( 2 {z}^{2}{\alpha}^{
	2}-2 {z}^{2}\alpha+{z}^{2}-1+W  ) $.

Plugging these results in and simplifying, we find explicit expressions for all four generating functions of interest, again with a common denominator $M$:

$
M= ( -1+z  )   ( z+1  )   ( -1+\alpha  ) ^{3}
 ( -1+{z}^{2}-3 {z}^{2}\alpha+3 {z}^{2}{\alpha}^{2}  )   ( {u}^{
	2}-{z}^{2}{u}^{4}\alpha-{z}^{2}{u}^{2}+2 {u}^{2}\alpha{z}^{2}+{\alpha}^{2}{z}^{2}{u}
^{4}-2 {u}^{2}{\alpha}^{2}{z}^{2}-{z}^{2}\alpha+{z}^{2}{\alpha}^{2}  ) 
$.

The first function:

$
4z^6\alpha^3MF_e=-  ( 3 {z}^{4}{\alpha}^{2}-{\alpha}^{2}W{z}^{2}-3 {z}^{2}{\alpha}^{2}-3 {z}^{4
}\alpha+aW{z}^{2}+3 {z}^{2}\alpha+{z}^{4}-W{z}^{2}-2 {z}^{2}+W+1  ) 
 ( -{z}^{4}{\alpha}^{4}-{z}^{4}{\alpha}^{2}+2 {z}^{4}{\alpha}^{3}-23 {z}^{4}{u
}^{2}{\alpha}^{2}+12 {z}^{4}{u}^{2}\alpha+22 {z}^{4}{u}^{2}{\alpha}^{3}-14 {\alpha}^{4}
{z}^{4}{u}^{2}-9 {z}^{4}{u}^{4}{\alpha}^{3}-2 {z}^{4}{u}^{4}\alpha+6 {z}^{4}{
	u}^{4}{\alpha}^{2}+5 {z}^{4}{u}^{4}{\alpha}^{4}-3 {z}^{4}{u}^{2}-6 {z}^{6}{\alpha}
^{5}-4 {z}^{6}{\alpha}^{3}+7 {z}^{6}{\alpha}^{4}+2 {\alpha}^{6}{z}^{6}+{z}^{6}{\alpha}^
{2}-6 {u}^{2}{z}^{6}{\alpha}^{5}-22 {u}^{2}{z}^{6}{\alpha}^{3}+16 {u}^{2}{z}^
{6}{\alpha}^{4}+5 {u}^{2}{\alpha}^{6}{z}^{6}+16 {u}^{2}{z}^{6}{\alpha}^{2}-6 {u}^{
	2}{z}^{6}\alpha-5 {\alpha}^{6}{z}^{6}{u}^{4}-13 {\alpha}^{4}{z}^{6}{u}^{4}+11 {\alpha}^
{5}{z}^{6}{u}^{4}+11 {\alpha}^{3}{z}^{6}{u}^{4}-5 {\alpha}^{2}{z}^{6}{u}^{4}+2
 {\alpha}^{7}{z}^{8}{u}^{4}+{\alpha}^{5}{z}^{8}{u}^{4}-3 {\alpha}^{6}{z}^{8}{u}^{4}
+{z}^{6}{u}^{4}\alpha+{u}^{2}{z}^{6}-2 {u}^{2}{\alpha}^{7}{z}^{8}+{u}^{2}{\alpha}^{6
}{z}^{8}+{z}^{2}{u}^{4}\alpha+3 {z}^{2}{u}^{2}-{\alpha}^{2}{z}^{2}{u}^{4}-6 {u
}^{2}\alpha{z}^{2}+7 {u}^{2}{\alpha}^{2}{z}^{2}-{u}^{2}+4 {z}^{2}{u}^{2}aW-{z}
^{2}{u}^{4}aW+{z}^{2}{u}^{4}{\alpha}^{2}W-5 {u}^{2}{\alpha}^{2}{z}^{2}W-8 {u}^
{2}{\alpha}^{3}W{z}^{4}+6 {u}^{2}{\alpha}^{4}W{z}^{4}+8 {u}^{2}{\alpha}^{2}W{z}^{4}
-4 {u}^{2}W{z}^{4}\alpha-3 {\alpha}^{2}{z}^{4}{u}^{4}W+5 {\alpha}^{3}{z}^{4}{u}^{4
}W-3 {\alpha}^{4}{z}^{4}{u}^{4}W+{\alpha}^{6}{z}^{6}{u}^{4}W+{z}^{4}{u}^{4}aW-{
	z}^{6}{u}^{4}{\alpha}^{5}W-{u}^{2}{\alpha}^{6}{z}^{6}W+{u}^{2}{z}^{4}W+{\alpha}^{2}W{
	z}^{4}-2 {\alpha}^{3}W{z}^{4}-2 {z}^{2}{u}^{2}W+{\alpha}^{4}W{z}^{4}+{u}^{2}W
 ) 
$.

The second function:

$
4z^5\alpha^3/(\alpha-1)/uMF_o= ( -3 {z}^{4}{\alpha}^{2}+{\alpha}^{2}W{z}^{2}+3 {z}^{2}{\alpha}^{2}+3 {z}^{4
}\alpha-aW{z}^{2}-3 {z}^{2}\alpha-{z}^{4}+W{z}^{2}+2 {z}^{2}-1-W  ) 
 ( 1+12 {z}^{4}{\alpha}^{4}+27 {z}^{4}{\alpha}^{2}-25 {z}^{4}{\alpha}^{3}-14
 {z}^{4}\alpha-{z}^{6}-7 {z}^{4}{u}^{2}{\alpha}^{2}+2 {z}^{4}{u}^{2}\alpha+10 {z}
^{4}{u}^{2}{\alpha}^{3}-5 {\alpha}^{4}{z}^{4}{u}^{2}+3 {z}^{4}+8 {z}^{6}{\alpha}^{
	5}+29 {z}^{6}{\alpha}^{3}-22 {z}^{6}{\alpha}^{4}-2 {\alpha}^{6}{z}^{6}-20 {z}^{6}
{\alpha}^{2}+7 {z}^{6}\alpha+5 W{z}^{4}\alpha-{z}^{4}W-5 aW{z}^{2}+5 {\alpha}^{2}W{z}^
{2}+2 W{z}^{2}-W-8 {u}^{2}{z}^{6}{\alpha}^{5}-14 {u}^{2}{z}^{6}{\alpha}^{3}+
15 {u}^{2}{z}^{6}{\alpha}^{4}+2 {u}^{2}{\alpha}^{6}{z}^{6}+6 {u}^{2}{z}^{6}{\alpha
}^{2}-{u}^{2}{z}^{6}\alpha+7 {z}^{2}\alpha-7 {z}^{2}{\alpha}^{2}-{u}^{2}\alpha{z}^{2}+{u
}^{2}{\alpha}^{2}{z}^{2}+{z}^{2}{u}^{2}aW-{u}^{2}{\alpha}^{2}{z}^{2}W-6 {u}^{2}
{\alpha}^{3}W{z}^{4}+3 {u}^{2}{\alpha}^{4}W{z}^{4}+4 {u}^{2}{\alpha}^{2}W{z}^{4}-{u
}^{2}W{z}^{4}\alpha-10 {\alpha}^{2}W{z}^{4}+9 {\alpha}^{3}W{z}^{4}-4 {\alpha}^{4}W{z}^{
	4}-3 {z}^{2}  )
$.

The third function:

$
8z^5\alpha^2MG_e= ( -3 {z}^{4}{\alpha}^{2}+{\alpha}^{2}W{z}^{2}+3 {z}^{2}{\alpha}^{2}+3 {z}^{4
}\alpha-aW{z}^{2}-3 {z}^{2}\alpha-{z}^{4}+W{z}^{2}+2 {z}^{2}-1-W  ) 
 ( 2 {\alpha}^{4}{z}^{4}{u}^{2}-2 {\alpha}^{4}{u}^{4}{z}^{4}-6 {z}^{4}{u
}^{2}{\alpha}^{3}+2 {z}^{4}{\alpha}^{3}+4 {\alpha}^{3}{u}^{4}{z}^{4}-{z}^{4}{\alpha}^{2}
-2 {\alpha}^{2}{u}^{4}{z}^{4}+{\alpha}^{2}W{z}^{2}+6 {z}^{4}{u}^{2}{\alpha}^{2}-3 
{z}^{2}{\alpha}^{2}-2 {z}^{4}{u}^{2}\alpha+2 {z}^{2}\alpha+1-W-{z}^{2}  ) 
 ( 2 {z}^{2}{\alpha}^{2}-2 {z}^{2}\alpha+{z}^{2}-1+W  )
$.

The fourth function:

$
8z^6(1-\alpha)\alpha^2/uMG_0= ( 2 {z}^{2}{\alpha}^{2}-2 {z}^{2}\alpha+{z}^{2}-1+W  )   ( 2 {
	z}^{4}{u}^{2}{\alpha}^{3}-2 {z}^{4}{\alpha}^{3}-3 {z}^{4}{u}^{2}{\alpha}^{2}+{z}^{4
}{\alpha}^{2}-{\alpha}^{2}W{z}^{2}+3 {z}^{2}{\alpha}^{2}-{u}^{2}{\alpha}^{2}{z}^{2}+{u}^{
	2}{\alpha}^{2}{z}^{2}W+{z}^{4}{u}^{2}\alpha+{u}^{2}\alpha{z}^{2}-{z}^{2}{u}^{2}\alpha W-2 
{z}^{2}\alpha+{z}^{2}+W-1  )   ( z\alpha+1  )   ( -1+z\alpha
 )   ( -3 {z}^{4}{\alpha}^{2}+{\alpha}^{2}W{z}^{2}+3 {z}^{2}{\alpha}^{2}+
3 {z}^{4}\alpha-\alpha W{z}^{2}-3 {z}^{2}\alpha-{z}^{4}+W{z}^{2}+2 {z}^{2}-1-W
 ) 
$.

Of course, the expressions do not look appealing, but that is what they are. We can derive as many corollaries 
from this as we want, of course with Maple:
\begin{align*}
f_0&=[u^0]F_e=1+ ( 2{\alpha}^{2}+1-2\alpha ) {z}^{2}+ ( 5{
	\alpha}^{4}-10{\alpha}^{3}+9{\alpha}^{2}-4\alpha+1 ) {z}^{
	4}+\cdots,\\
f_1&=[u^1]F_o=\alpha z+ ( 3{\alpha}^{2}-4\alpha+2 ) \alpha{z}^{3}
+ ( 8{\alpha}^{4}-19{\alpha}^{3}+20{\alpha}^{2}-11\alpha+
3 ) \alpha{z}^{5}+\cdots,\\
f_2&=[u^2]F_e=\alpha ( 1-\alpha ) {z}^{2}+2 ( 1-\alpha
)  ( 2{\alpha}^{2}+1-2\alpha ) \alpha{z}^{4}+\cdots,\\
f_3&=[u^3]F_o= ( 1-\alpha ) {\alpha}^{2}{z}^{3}+ (1-\alpha
)  ( 5{\alpha}^{2}-6\alpha+3 ) {\alpha}^{2}{z}^{
	5}+\cdots.
\end{align*}
and similarly
\begin{align*}
	g_0&=[u^0]G_e=\alpha ( 1-\alpha ) ^{2}{z}^{3}+ ( 5{\alpha}^{2}
	-4\alpha+2 )  ( 1-\alpha ) ^{2}\alpha{z}^{5}
	+\cdots,\\
	g_1&=[u^1]G_o=\alpha ( 1-\alpha ) {z}^{2}+2 ( 1-\alpha
	)  ( 2{\alpha}^{2}+1-2\alpha ) \alpha{z}^{4}
	,\\
	g_2&=[u^2]G_e= ( 1-\alpha ) {\alpha}^{2}{z}^{3}+ ( 1-\alpha
	)  ( 5{\alpha}^{2}-6\alpha+3 ) {\alpha}^{2}{z}^{		5}
	+\cdots,\\
	g_3&=[u^3]G_o={\alpha}^{2} ( 1-\alpha ) ^{2}{z}^{4}+3{\alpha}^{2}
	( 2{\alpha}^{2}+1-2\alpha )  ( 1-\alpha ) 	^{2}{z}^{6}
	+\cdots.
\end{align*}

\section{A more sophisticated approach}

The imbalance of $\alpha$ versus $\beta$ is leveled out after 2 (or an even number of) steps. Thus, as in \cite{PP}, we consider the system after an even number of steps.
In the following graph, a directed arrow stands for 2 steps (a double-step). Note that the system is still working without look-ahead, writing  \textsf{s}
for the small item of size $\frac13$ and \textsf{l}
for the large item of size $\frac23$, the sequences \textsf{sl} resp.\ \textsf{ls} lead to different states when being in the special state named $Q$.

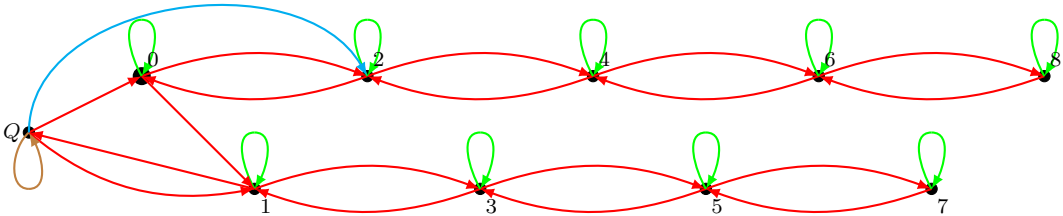
\begin{figure}[h]

	\begin{center}
		\begin{tikzpicture}[scale=1.5]
			%\draw (0,0) circle (0.1cm);

			\foreach \x in {0,2,4,6,8}
			{
				\draw (\x,0) circle (0.05cm);
				\fill (\x,0) circle (0.05cm);
				%	\node[circle,draw ] (c) at (\x,0){$\x$};
			}
			
			\foreach \x in {1,3,5,7}
			{
				\draw (\x,-1) circle (0.05cm);
				\fill (\x,-1) circle (0.05cm);
				%	\node[circle,draw ] (c) at (\x,0){$\x$};
			}

			\draw (-1,-0.5) circle (0.05cm);	\fill (-1,-0.5) circle (0.05cm);

%			\draw[thick, cyan, -latex] (0,0) to [out=150,in=-30] (-1,0.5);	
	%		\draw[thick, cyan, -latex]  (-1,0.5)to [out=00,in=120](0,0);	

		%	\draw[thick, red, -latex] (0,-1) to  (-1,-1.5);	
		%	\draw[thick, red, -latex]  (-1,-1.5)to [out=0,in=-120](0,-1);	

			%\draw[thick, red, -latex]    (0,-1) to (-1,0.5);
		%	\draw[thick, cyan, -latex]     (-1,-1.5)to (1,0) ;
		%	\draw[thick, red, -latex]     (-1,0.5)to (1,-1) ;
			
			\fill (0,0) circle (0.08cm);

			\foreach \x in {0,2,4,6}
			{
				\draw[thick, red, -latex] (\x,0) to[out=20,in=160]  (\x+2,0);	
				\draw[thick, red, -latex] (\x+2,0) to[out=200,in=-20]  (\x,0);	
				\draw[thick, green, ] (\x,0) to[out=120,in=180]  (\x,0.5);	
				\draw[thick, green, -latex] (\x,0.5) to[out=0,in=60]  (\x,0);	
			}
		\draw[thick, green, ] (8,0) to[out=120,in=180]  (8,0.5);	
		\draw[thick, green, -latex] (8,0.5) to[out=0,in=60]  (8,0);	
			%\foreach \x in {1,3,5,7}
			%{
		%		\draw[thick, cyan, -latex] (\x,0) to[out=20,in=160]  (\x+1,0);	
		%		\draw[thick, cyan, -latex] (\x+1,0) to[out=200,in=-20]  (\x,0);	
		%	}

			%\foreach \x in {0,2,4,6}
			%{
		%		\draw[thick, cyan, -latex] (\x,-1) to[out=20,in=160]  (\x+1,-1);	
		%		\draw[thick, cyan, -latex] (\x+1,-1) to[out=200,in=-20]  (\x,-1);	
	%		}
			
			\foreach \x in {1,3,5}
			{
				\draw[thick, red, -latex] (\x,-1) to[out=20,in=160]  (\x+2,-1);	
				\draw[thick, red, -latex] (\x+2,-1) to[out=200,in=-20]  (\x,-1);	
				\draw[thick, green, ] (\x,-1) to[out=120,in=180]  (\x,-0.5);	
				\draw[thick, green, -latex] (\x,-0.5) to[out=0,in=60]  (\x,-1);
			}
				\draw[thick, green, ] (7,-1) to[out=120,in=180]  (7,-0.5);	
			\draw[thick, green, -latex] (7,-0.5) to[out=0,in=60]  (7,-1);
			
			\foreach \x in {0,2,4,6}
			{
				%	\draw[thick, -latex] (\x+1,0) to  (\x,0);	
				\node at  (\x+0.1,0.15){\tiny$\x$};
			}			
			\foreach \x in {1,3,5,7}
			{
				%	\draw[thick, -latex] (\x+1,0) to  (\x,0);	
				\node at  (\x+0.1,-1.15){\tiny$\x$};
			}
			
			\draw[thick, red, -latex] (0,0) to  (1,-1);
				\draw[thick, red, -latex] (1,-1) to  (-1,-0.5);
				
					\draw[thick, red, -latex] (-1,-0.5) to[out=-40,in=190] (1,-1);
			
			\draw[thick, brown ] (-1,-0.5) to[out=220,in=180]  (-1,-1);	
			\draw[thick, brown, -latex] (-1,-1) to[out=0,in=-60]  (-1,-0.5);
			\draw[thick, red,-latex ] (-1,-0.5) to  (0,0);	
						\draw[thick, cyan,-latex ] (-1,-0.5) to[out=90,in=120]  (2,0);	
			
			\node at  (8+0.1,0.15){\tiny$8$};
			\node at  (-1.15,-0.5){\tiny$Q$};
			
		\end{tikzpicture}
	\end{center}
	\caption{Two steps. Red with probability $\alpha\beta$, green with probability $1-2\alpha\beta=\alpha^2+\beta^2$, blue with probability $\beta^2$, brown with probability $\alpha^2$.}
\end{figure}

The graph is now simpler than before. We introduce generating functions $f_N$ for the upper layer of states, and $g_N$ for the lower layer of states.
The meaning of these generating functions is now different from the previous section, but it is apparent how they are related. Here are the recursions:
\begin{align*}
f_N&=z\alpha\beta f_{N-1}+z\alpha\beta f_{N+1}+z(\alpha^2+\beta^2)f_N,\quad N\ge2,\\
f_1&=z\alpha\beta f_{0}+z\alpha\beta f_{2}+z(\alpha^2+\beta^2)f_1+z\beta^2f_Q,\\
f_0&=1+z\alpha\beta f_{1}+z(\alpha^2+\beta^2)f_0+z\alpha\beta f_Q,\\
f_Q&=z\alpha\beta g_0+z\alpha^2f_Q=\frac{z\alpha\beta g_0}{1-z\alpha^2},\\
g_N&=z\alpha\beta g_{N-1}+z\alpha\beta g_{N+1}+z(\alpha^2+\beta^2)g_N,\quad N\ge1,\\
g_0&=z\alpha\beta f_{0}+z\alpha\beta g_{1}+z\alpha\beta f_{Q}+z(\alpha^2+\beta^2)g_0.
\end{align*}
Introducing only \emph{two} bivariate generating functions 
\begin{equation*}
F(u)=\sum_{N\ge0}u^Nf_N\quad\text{and}\quad G(u)=\sum_{N\ge0}u^Ng_N, 
\end{equation*}
we find by summing the recursions
\begin{align*}
F(u)=\sum_{N\ge0}u^Nf_N&=z\alpha\beta \sum_{N\ge2}u^Nf_{N-1}+z\alpha\beta \sum_{N\ge2}u^Nf_{N+1}+z(\alpha^2+\beta^2)\sum_{N\ge2}u^Nf_N
\\&+u(z\alpha\beta f_{0}+z\alpha\beta f_{2}+z(\alpha^2+\beta^2)f_1+z\beta^2\frac{z\alpha\beta g_0}{1-z\alpha^2})\\
&+1+z\alpha\beta f_{1}+z(\alpha^2+\beta^2)f_0+z\alpha\beta \frac{z\alpha\beta g_0}{1-z\alpha^2}\\
&=z\alpha\beta uF(u) +\frac{z\alpha\beta}{u} (F(u)-f_0)+z(\alpha^2+\beta^2)F(u)
\\&+uz^2\beta^3\alpha\frac{ g_0}{1-z\alpha^2}+1+z^2\alpha^2\beta^2 \frac{ g_0}{1-z\alpha^2}\\
\end{align*}
and
\begin{align*}
	G(u)=\sum_{N\ge0}u^Ng_N&=z\alpha\beta \sum_{N\ge1}u^Ng_{N-1}+z\alpha\beta\sum_{N\ge1}u^N g_{N+1}+z(\alpha^2+\beta^2)\sum_{N\ge1}u^Ng_N\\
	&+z\alpha\beta f_{0}+z\alpha\beta g_{1}+z\alpha\beta \frac{z\alpha\beta g_0}{1-z\alpha^2}+z(\alpha^2+\beta^2)g_0\\
	&=z\alpha\beta u G(u)+\frac{z\alpha\beta}{u}(G(u)-g_0)+z(\alpha^2+\beta^2)G(u)\\
	&+z\alpha\beta f_{0}+z^2\alpha^2\beta^2\frac{ g_0}{1-z\alpha^2}.
\end{align*}
Solving the system leads to
\begin{align*}
F(u)&={\frac {-uz{\alpha}^{2}+{z}^{2}{\alpha}^{2}{\beta}^{2}g_0 u-z
		\alpha \beta f_0+{z}^{2}{\alpha}^{3}\beta f_0+u+{u}^{2}{z
		}^{2}{\beta}^{3}\alpha g_0}{u-2 uz{\alpha}^{2}-z\alpha \beta 
		{u}^{2}+{z}^{2}{\alpha}^{3}\beta {u}^{2}-z\alpha \beta+{z}^{2}{
			\alpha}^{3}\beta+{z}^{2}u{\alpha}^{4}-zu{\beta}^{2}+{z}^{2}u{\beta}^{2
		}{\alpha}^{2}}},\\
	G(u)&=-{\frac {z\alpha \beta   ( -z\alpha \beta g_0 u+g_0
			-g_0 z{\alpha}^{2}+z{\alpha}^{2}f_0 u-f_0 u  ) 
		}{u-2 uz{\alpha}^{2}-z\alpha \beta {u}^{2}+{z}^{2}{\alpha}^{3}\beta
			 {u}^{2}-z\alpha \beta+{z}^{2}{\alpha}^{3}\beta+{z}^{2}u{\alpha}^{4}
			-zu{\beta}^{2}+{z}^{2}u{\beta}^{2}{\alpha}^{2}}}.
	\end{align*}
These answers are implicit, since they contain $f_0=F(0)$ and $g_0=G(0)$. To make them explicit, the kernel method is used once again. 
The denominators factor as
\begin{equation*}
z\alpha \beta   ( -1+z{\alpha}^{2}  )(u-r_1)(u-r_2)
\end{equation*}
with
\begin{align*}
r_2={\frac {1-z{\alpha}^{2}-z{\beta}^{2}-\sqrt {{z}^{2}{\alpha}^{4}-
			2 {z}^{2}{\beta}^{2}{\alpha}^{2}-2 z{\alpha}^{2}+{z}^{2}{\beta}^{4}-
			2 z{\beta}^{2}+1}}{2z\alpha \beta}}
\end{align*}
and $r_1=\frac1{r_2}$.

The factor $(u-r_2)$ (the `bad' factor) must cancel from numerator and denominator. The result is now
\begin{equation*}
F(u)=\frac{r_2 {z}^{2}{\beta}^{3}\alpha g_0-z{\alpha}^{2}+{z}^{2}{
		\alpha}^{2}{\beta}^{2}g_0+1+u{z}^{2}{\beta}^{3}\alpha g_0
}{z\alpha \beta   ( -1+z{\alpha}^{2}  )(u-r_1)}
\end{equation*}
and
\begin{equation*}
	G(u)=\frac{-  ( -z\alpha \beta g_0+z{\alpha}^{2}f_0-f_0
		 ) 		}{  ( -1+z{\alpha}^{2}  )(u-r_1)},
\end{equation*}
Plugging in $u=0$, we get
\begin{align*}
f_0&=\frac{r_2 {z}^{2}{\beta}^{3}\alpha g_0-z{\alpha}^{2}+{z}^{2}{		\alpha}^{2}{\beta}^{2}g_0+1
}{z\alpha \beta   ( -1+z{\alpha}^{2}  )(-r_1)},\\
g_0&=\frac{  ( -z\alpha \beta g_0+z{\alpha}^{2}f_0-f_0
	 ) 
}{ ( -1+z{\alpha}^{2}  )r_1}.
\end{align*}
From these, we can compute $f_0$ and $g_0$ easily, but don't print it, since it is not too attractive at the moment (in a moment, it will become very beautiful).

It is easy to see that
\begin{equation*}
[u^j]G(u)=\frac{  ( -z\alpha \beta g_0+z{\alpha}^{2}f_0-f_0
	 ) 		}{  ( -1+z{\alpha}^{2}  )}r_2^{j+1}
\end{equation*}
and
\begin{align*}
[u^j]F(u)&=-r_2^{j+1}\frac{r_2 {z}^{2}{\beta}^{3}\alpha g_0-z{\alpha}^{2}+{z}^{2}{
		\alpha}^{2}{\beta}^{2}g_0+1
}{z\alpha \beta   ( -1+z{\alpha}^{2}  )}
-r_2^{j}\frac{{z}{\beta}^{2}g_0
}{  ( -1+z{\alpha}^{2}  )}.
\end{align*}
Note that $[z^mu^j]F(u)$ is the probability to reach state $2j$ in $m$ (double-)steps, and
$[z^mu^j]G(u)$ is the probability to reach state $2j+1$ in $m$ (double-)steps.

\subsection*{More attractive formul\ae\ thanks to a substitution}

Using the substitution
\begin{equation*}
z=\frac{v}{\alpha\beta+(\alpha^2+\beta^2)v+\alpha\beta v^2}=\frac{v}{  ( \alpha+v\beta  )( \beta+v\alpha  )  },
\end{equation*}
(inspired by our old paper \cite{PP})
all the expressions become nicer.\footnote{A referee pointed out the similarity to the \emph{Joukowsky transform}, https://en.wikipedia.org/wiki/Joukowsky_transform.} For instance, $r_2=v$ and
\begin{align*}
f_0&={\frac {  ( v\alpha+\beta  )   ( \alpha+v\beta  ) }
	{\alpha \beta  (1-v)   ( {v}^{2}+v+1  ) }},\\
g_0&={\frac {v  ( \alpha+\alpha {v}^{2}+v\beta  )   ( v
		\alpha+\beta  ) }{\alpha \beta   (1-v)   ( {v}^		{2}+v+1  ) }}.
\end{align*}
The equality $(1-v)   ( {v}^		{2}+v+1  )=1-v^3$ might be useful as well.
Even the full bivariate generating functions look now very nice:
\begin{align*}
F&={\frac {  ( u{v}^{3}\beta+\alpha+v\beta  )   ( v\alpha+
		\beta  ) }{\beta \alpha   ( 1-uv  )   ( 1-v		 )   ( {v}^{2}+v+1  ) }},\\
G&={\frac {v  ( \alpha+\alpha {v}^{2}+v\beta  )   ( v\alpha
		+\beta  ) }{\beta \alpha   ( 1-uv  )   ( 1-v		 )   ( {v}^{2}+v+1  ) }}.
\end{align*}
Consequently, reading off coefficient of powers of $u$,
\begin{equation*}
[u^j]	F={\frac {v^j  (\alpha+v\beta  )   ( v\alpha+
			\beta  ) }{\beta \alpha     ( 1-v		 )   ( {v}^{2}+v+1  ) }}
		+{\frac {v^{j+1} }{\alpha     ( 1-v		 )   ( {v}^{2}+v+1  ) }}
\end{equation*}
and
\begin{equation*}
[u^j]	G={\frac {v^{j+1}  ( \alpha+\alpha {v}^{2}+v\beta  )   ( v\alpha
			+\beta  ) }{\beta \alpha     ( 1-v		 )   ( {v}^{2}+v+1  ) }}.
\end{equation*}
Finally we answer the question how to read off coefficients of powers of $z$ when the function is given in terms of $v$: For that, we employ Cauchy's integral formula in the following computation,
\begin{align*}
[z^N]H(z(v))&=\frac1{2\pi i}\oint \frac{dz}{z^{N+1}}H(z(v))\\
&=\frac1{2\pi i}\oint \frac{dv}{v^{N+1}}\frac{\alpha\beta(1-v^2)}{(\alpha+\beta v)(\beta+\alpha v)}(\alpha+\beta v)^{N+1}(\beta+\alpha v)^{N+1}H(v)\\
&=[v^N]\alpha\beta(1-v^2)(\alpha+\beta v)^{N}(\beta+\alpha v)^{N}H(v).
\end{align*}

\subsection*{Walks with an odd number of steps}

For that, we don't need to do new calculations, by considering the last step separately. We refer to the original Figure~\ref{KBH}. It is immediate to see that
\begin{align*}
\mathbb{P}\{\text{reach top level state } &2j+1 \text{ in $2m+1$ steps} \}\\&=
\alpha \mathbb{P}\{\text{reach top level state $2j$ in $2m$ steps} \}\\
&+\beta \mathbb{P}\{\text{reach top level state $2j+2$ in $2m$ steps} \}
\end{align*}
and
\begin{align*}
	\mathbb{P}\{\text{reach bottom  level state } &2j \text{ in $2m+1$ steps} \}\\&=
	\alpha \mathbb{P}\{\text{reach bottom level state $2j-1$ in $2m$ steps} \}\\
	&+\beta \mathbb{P}\{\text{reach bottom level state $2j+1$ in $2m$ steps} \};
\end{align*}
the exceptional cases near the beginning are easy to figure out directly.

\section{Asymptotics}

Although this paper concentrates on \emph{explicit enumerations}, one referee asks for some asymptotic
considerations. One natural concept would the height of a Kn\"odel walk, i.~e., the state with the highest index that is
reached during the walk. For simpler walks, this has been worked out in \cite{Knoedel-height}, compare \cite{PP}.
However, that would be a completely different approach, and we have chosen the kernel method as the unifying 
method of choice.

We can, however, offer something appealing here, namely we compute the average index of the state where the walks ends
(in the sophisticated version). So we compute
\begin{equation*}
\textsc{expected-end}=\sum_{k\ge0}(2k) f_k+\sum_{k\ge0}(2k+1) g_k.
\end{equation*}
This is best computed using the bivarite generating functions:
\begin{align*}
	\textsc{expected-end}&=\frac{\partial}{\partial u}\Big(F(u^2,z)+uG(u^2,z)\Big)\bigg|_{u=1}\\
	&={\frac {v \left( v\alpha+\beta \right)  \left( 3{v}^{2}\beta+3
			\alpha+3v\beta+v\alpha+\alpha{v}^{2}+\alpha{v}^{3} \right) }{
		\alpha 	\beta\left( 1-v \right) ^{3} \left( 1+v+{v}^{2} \right) }}.
	\end{align*}
To find asymptotics, we transfer back from $v$ to $z$. In order to avoid ungainly expressions that the reader can generate himself/herself with a computer, 
we demonstrate the procedure for the standard case $\alpha=\beta=\frac12$. Then 
$$v=\dfrac{-z+2-2\sqrt{1-z}}{z}\sim \frac{z}{4}+\frac{z^2}{8}+\cdots$$ and
\begin{align*}
	\textsc{expected-end}&=\frac{z-1+(1+z)\sqrt{1-z}}{2(1-z)^2}\sim\frac{1}{(1-z)^{3/2}}.
	\end{align*}
The coefficient of $z^n$ in this expression is $\binom{-3/2}{n}(-1)^n\sim 2\sqrt{\frac n\pi}$, which is the answer to the question about the average index where
the walk stops. In the general instance, the singularity of interest is $z\sim1$, which is equivalent to $v\sim1$, and 
\begin{equation*}
1-v\sim \frac{1}{\sqrt{\alpha\beta}}\sqrt{1-z}.
\end{equation*}
But
\begin{align*}
	\textsc{expected-end} \sim\frac2{\alpha\beta}\frac1{(1-v)^3}\sim2\sqrt{\alpha\beta} \frac{1}{(1-z)^{3/2}}, 
\end{align*}
whence the  result in the general case is $\sim4\sqrt{\alpha\beta}\sqrt{\frac n\pi}$.

The method of local expansions around the dominant singularity (here $z=1$) and then translating into the behaviour of the coefficients is called
\emph{singularity analysis of generating functions}. Standard references are \cite{FS} and \cite{FO}.

\section{Conclusion}

We want to emphasize the following points:

\begin{itemize}
	\item A brute-force approach is possible, but leads to equations of order 4 and explicit but very ungainly expressions.
	
	\item Looking at the system after an even number of steps is a clever idea, since the imbalance of $\alpha$ versus $\beta$ is leveled out. The equations are only quadratic.
	
	\item Introducing an auxiliary variable, all the generating functions become rational (in the variable). Consequently reading off coefficients is not difficult.
	
	\item To go from an even number of steps to an odd number of steps is not difficult, when considering the last step separately and use previous results.
	
	\item Once the generating functions of interest are known explicitly, several corollaries of an asymptotic nature can be derived from them. 
\end{itemize}

\bibliographystyle{plain}

%\bibliography{hornik2}

\end{document}